% Bivariant Chern classes and Grothendieck transformations
\input amstex
\documentstyle{amsppt}
\pagewidth{13.5cm}
\magnification=1200
\pageheight{19cm}
%\input xypic
%%%%%%%%%%%%%%%%%%%%%%%%%%%%%%%%%%%%%%%%%%%%%%%%
\define\alp{\alpha}
\define\be{\beta}
\def \jeden {1\hskip-3.5pt1}

\define\ga{\gamma}

\define\de{\delta}

\define\bA{\Bbb A}
\define\bF{\Bbb F}
\define\bH{\Bbb H}
\define\bPH{\Bbb {PH}}
\define\bSH{\Bbb {SH}}
\define\bK{\Bbb K}
\define \bB{\Bbb B}

\define \bZ{\Bbb Z}
\define \bC{\Bbb C}

\define \bQ{\Bbb Q}
\def\op{\operatorname}

%%%%%%%%%%%%%%%%%%%%%%%%%%%%%%%%%%%%%%%%%%%%%%%%
\topmatter
\title Bivariant Chern classes and Grothendieck transformations
\endtitle
\author Jean-Paul Brasselet, J\"org Sch\"urmann and Shoji Yokura$^{(*)}$
\endauthor
\thanks {(*) Partially supported by Grant-in-Aid for Scientific Research (C) (No. 15540086), the Japanese Ministry of Education, Science, Sports and Culture} 
\endthanks
\address
Institute de Math\'ematiques de Luminy, UPR 9016 - CNRS, Campus de Luminy - Case 907, 13288 Marseille Cedex 9, France
\endaddress
\email jpb$\@$iml.univ-mrs.fr
\endemail
\address
Westf. Wilhelm-Universit\"at, SFB 478 ``Geometrische Strukturen in der Mathematik", Hittorfstr. 27, 48149 M\"unster, Germany
\endaddress
\email jschuerm$\@$math.uni-muenster.de
\endemail
\address
Department of Mathematics and Computer Science, 
Faculty of Science, 
University of Kagoshima, 21-35 Korimoto 1-chome, Kagoshima 890-0065, Japan
\endaddress
\email yokura$\@$sci.kagoshima-u.ac.jp
\endemail
\abstract {The existence of bivariant Chern classes was conjectured by W. Fulton and R. MacPherson and 
proved by J.-P. Brasselet for cellular morphisms of analytic varieties. 
However, its uniqueness has been unresolved since then. In this paper we show that restricted to morphisms 
whose target varieties are possibly singular but (rational) homology manifolds (such as orbifolds), 
the bivariant Chern classes (with rational coefficients) are uniquely determined. 
And also we discuss some related things and problems. In the final sections we construct a unique bivariant 
Chern class $\widetilde \ga : \widetilde {\bF} \to \widetilde {\bH}$ satisfying a suitable normalization condition. In fact, it will be a special case of a general construction of unique Grothendieck transformations, which in a sense gives a positive answer to the uniqueness questions concerning Grothendieck transformations posed by Fulton and MacPherson.
} 
\endabstract
\endtopmatter

\document
\head \S 1 Introduction \endhead

Various characteristic classes of singular varieties
have been introduced and studied. One of them is the so-called 
Chern--Schwartz--MacPherson class [Mac], which seems to be a fundamental and important 
characteristic class. Its unique existence was conjectured by P. Deligne 
and A. Grothendieck and it was finally solved affirmatively by R. 
MacPherson [Mac]. It turns out that the Chern--Schwartz--MacPherson class is 
isomorphic to the Schwartz class ([Schw1, Schw2]) by the Alexander duality isomorphism [BS], 
which is a reason for the name (at least for compact complex varieties embeddable into manifolds).

In early 1980, W. Fulton and R. MacPherson introduced the notion of 
bivariant theory which associates to a morphism an abelian group, unifying 
covariant and contravariant theories. In particular, they asked if there 
exists a (unique) Grothendieck transformation $\ga: \bF \to \bH$ from the 
bivariant theory $\bF$ of constructible functions to 
the bivariant homology theory $\bH$ such that it specializes to the 
original Chern--Schwartz--MacPherson class for a morphism to a point.
Such a transformation is called a {\it bivariant Chern class} for short.

In [B] J.-P. Brasselet shows that a bivariant Chern class exists for cellular 
morphisms between complex varieties embeddable into complex manifolds. Since then its uniqueness has been unresolved (cf. [Z1, Z2]). It is showed in [Y3] that 
in the case when the target variety is nonsingular, the bivariant Chern class is uniquely determined;
more precisely, if there exists a bivariant Chern class $\ga :\bF \to \bH$, 
then for a morphism $f: X \to Y$ with $Y$ being nonsingular and for 
a bivariant constructible function $\alp$ the following holds:
$$\ga(\alp) = f^*s(TY) \cap c_*(\alp)$$
with $s(TY) := c^*(TY)^{-1}$ the total Segre class of the tangent bundle $TY$.
Conversely, it is showed in [Y4] that when restricted to morphisms with nonsingular target varieties there 
exists a bivariant theory $\widetilde {\bF}$ of constructible functions (which conjecturally contains $\bF$) and
the transformation $\ga^{\op {Gin}}:= f^*s(TY) \cap c_* : \widetilde {\bF}(X \overset f \to \longrightarrow Y) \to \bH(X \overset f \to \longrightarrow Y)$ is the unique bivariant Chern class.

In this paper we show the uniqueness of bivariant Chern classes even when the target variety $Y$ is singular but a suitable homology manifold and discuss related problems. And in the final section we show that there exist a bivariant theory $\widetilde {\bF}$ of constructible functions, a bivariant homology theory $\widetilde {\bH}$ and a unique bivariant Chern class
$\widetilde \ga : \widetilde {\bF} \to \widetilde {\bH}$ satisfying the normalization condition. More precisely, we will show that a natural transformation of covariant theories extends uniquely to a Grothendieck transformation of suitable bivariant theories associated to the covariant theories, if the given transformation commutes with exterior products. This is in a sense a positive answer to [FM, \S10.9 Uniqueness questions].

\head \S2 Bivariant constructible functions and bivariant homology \endhead

For a general reference for the bivariant theory, see Fulton--MacPherson's book [FM]. 

For a category $\Cal C$ which has a final object $pt$ and on which the fiber product is well-defined, a bivariant theory $\bB$ on the category $\Cal C$ with values in the category of abelian groups is an assignment to each morphism
$$ X  \overset f \to \longrightarrow Y$$
in the category $\Cal C$ a (graded) abelian group
$$\bB(X  \overset f \to \longrightarrow Y)$$
which is equipped with the following three basic operations:

\noindent (Product operations): For morphisms $f: X \to Y$ and $g: Y
\to Z$, the product operation
$$\bullet: \bB( X \overset f \to \longrightarrow Y) \otimes \bB( Y \overset g
\to \longrightarrow Z) \to
\bB( X \overset {gf} \to \longrightarrow Z)$$
is  defined.

\noindent (Pushforward operations): For morphisms $f: X \to Y$
and $g: Y \to Z$ with $f$ confined, the pushforward operation
$$f_*: \bB( X \overset {gf} \to \longrightarrow Z) \to \bB( Y \overset g
\to \longrightarrow Z) $$
is  defined.

\noindent (Pullback operations): For an independent square (which we assume to be Cartesian)
$$\CD
X' @> g' >> X \\
@V f' VV @VV f V \\
Y' @> g >> Y, \endCD 
$$
the pullback operation
$$g^* : \bB( X \overset  f \to \longrightarrow Y) \to \bB( X' \overset g
\to \longrightarrow Y') $$
is  defined.

And these three operations are required to satisfy seven natural compatibility axioms [FM, Part I, \S 2.2].

$B_*(X) := \bB(X \to pt)$ becomes a covariant functor and $B^*(X) := \bB(X \overset {\op {id}} \to\longrightarrow X)$ becomes a contravariant functor.

Let $\bB, \bB'$ be two bivariant theories on two such categories $\Cal C$ and $\Cal C'$. Let $ - : \Cal C \to \Cal C'$ be a functor respecting the underlying structures. Then
a {\it Grothendieck transformation} from $\bB$ to $\bB'$
$$\ga : \bB \to \bB'$$
is a collection of homomorphisms
$$\ga_f: \bB(X \overset  f \to \longrightarrow Y) \to \bB'(\overline {X} \overset  \overline f \to \longrightarrow \overline Y) $$
for a morphism $X \overset  f \to \longrightarrow Y$ in the category $\Cal C$, which preserves the above three basic operations.

Note that using the target bivariant theory $\bB'$ one can define the abelian group
$$\bB''(X \overset  f \to \longrightarrow Y) := \bB'(\overline {X} \overset  \overline f \to \longrightarrow \overline Y).$$
It turns out that $\bB''$ becomes a bivariant theory on the category $\Cal C$ with the obvious bivariant operations. And the above Grothendieck transformation $\ga: \bB \to \bB'$ factors uniquely through
the tautological Grothendieck transformation $\ga^{\op {taut}}: \bB \to \bB''$. So, for the discussion of the
existence or uniqueness of a suitable Grothendieck transformation, it suffices to work out the case when
$ - : \Cal C \to \Cal C'$  is the identity $\op {id}: \Cal C \to \Cal C$, which we consider from now in the rest of the paper.
 
One important example for us is given by the theory of constructible functions. Let $X$ be a complex analytic (algebraic) set and let $F(X)$ denote the abelian group of all the complex analytically (algebraically) constructible functions on $X$. The association $X \longmapsto F(X)$ becomes a contravariant functor with the usual functional pullback. Moreover, for proper morphisms (which are confined maps in this context), it also becomes a covariant
functor with the pushforward $f_*$ defined by taking the weighted Euler--Poincar\'e 
characteristics fiberwise in the following sense:
$$f_*(\alpha)(y) = \chi \Bigl (f^{-1}(y);\alpha \Bigr ) \quad \text {for} \quad \alpha \in F(X) $$
$$\chi(A;\alpha) = \sum _{n \in \bZ} n \chi(A \cap \alp^{-1}(n)).$$
Furthermore this becomes a bivariant theory as follows:
For {\it any} morphism $f: X \to Y$ the group $s\bF(X \to Y)$ is defined by
$$s\bF(X \overset f \to \longrightarrow Y) := F(X).$$
and the three operations are defined as follows:

\noindent (i) the product operation
$$\bullet: s\bF( X \overset f \to \longrightarrow Y) \otimes s\bF( Y \overset g
\to \longrightarrow Z) \to
s\bF( X \overset {gf} \to \longrightarrow Z)$$
is  defined by:
$$\alp \bullet \be:= \alp \cdot f^*\be.$$

\noindent (ii) the pushforward operation for $f$ proper
$$f_*: s\bF( X \overset {gf} \to \longrightarrow Z) \to s\bF( Y \overset g
\to \longrightarrow Z) $$
is the pushforward
$$ f_* :F(X) \to F(Y).$$

\noindent (iii) For a fiber square (which is the independent square in this context)
$$\CD
X' @> g' >> X \\
@V f' VV @VV f V\\
Y' @> g >> Y, \endCD %\tag 2.0
$$
the pullback operation
$$g^* : s\bF( X \overset  f \to \longrightarrow Y) \to s\bF( X' \overset g
\to \longrightarrow Y') $$
is  the pullback 
$${g'}^*: F(X) \to F(X').$$
Then $s\bF$ becomes a bivariant theory and shall be called a {\it simple bivariant theory of constructible functions}.

Fulton--MacPherson's bivariant theory $\bF$ of constructible functions is much more subtle and requires some strong geometric or topological conditions coming from morphisms themselves. For a morphism $f: X \to Y$ the bivariant theory $\bF(X \overset f \to \longrightarrow Y)$ of constructible functions consists of all the constructible functions on $X$ which satisfy the {\it local Euler condition with respect to $f$}, i.e., the condition that for any point $x \in X$ and for any local embedding $(X, x) \to (\bold C^N,
0)$ the following equality holds
$$ \alpha(x) = \chi \bigl (B_{\epsilon} \cap f^{-1}(z); \alpha \bigr),$$
where $B_{\epsilon}$ is a sufficiently small open ball of the origin $0$
with radius $\epsilon$ and $z$ is any point close to $f(x)$ (see [B], [FM], [S]). 
The three operations on $\bF$ are the same as in $s\bF$, i.e., $\bF$ is a bivariant subtheory of $F$.
Note that $\bF(X \to pt) = F(X)$ and $\bF(X \overset {\op {id}_X} \to \longrightarrow X)$ consists of all locally constant functions. 

Let us recall briefly the definition of the Fulton--MacPherson bivariant homology theory $\bH$ constructed from the usual cohomology theory with values in a fixed commutative ring $A$ (see [FM, \S 3.1] ). For a morphism $f: X \to Y$, choose a morphism $\phi: X \to M$ to a smooth oriented manifold $M$ of real dimension $n$ such that 
$\Phi:= (f, \phi): X \to Y \times M$ is a closed embedding. Then the $i$-th bivariant homology group 
$\bH^i(X \overset f \to \longrightarrow Y)$ is defined by
$$  \bH^i(X \overset f \to \longrightarrow Y) := H^{i+n}(Y \times M, (Y \times M) \setminus X_{\phi}; A),$$
where $X_{\phi} = \Phi(X)$.  

A bivariant Chern class is a Grothendieck transformation from the bivariant theory $\bF$ of constructible functions 
to the bivariant homology theory $\bH$
$$\ga: \bF \to \bH$$
satisfying the normalization condition that
for a nonsingular variety $X$ and for the map $\pi: X \to pt$ to a point $pt$
$$\ga(\jeden_{\pi}) = c(TX) \cap [X] $$
where 
$\jeden_{\pi} = \jeden_X \in F(X) = \bF(X \overset {\pi} \to\longrightarrow pt)$, $c(TX) \in H^*(X)$ is the total Chern class of the tangent bundle $TX$ and $[X] \in H_*(X)$ is the fundamental class of $X$. 
For the existence of such a bivariant Chern class (under some mild restrictions), we recall the following theorem:
\proclaim {Theorem (2.1)} (Brasselet's Theorem [B]) For the category of embeddable analytic varieties with cellular morphisms there exists a bivariant Chern class $\ga: \bF \to \bH$.
\endproclaim
\head \S3 On the uniqueness of the bivariant Chern classes \endhead

First we recall the notion of strong orientation from [FM].

\definition {Definition (3.1)} An element $\theta \in \bB(X \overset f \to\longrightarrow Y)$ is called a strong 
orientaion for the morphism $f: X \to Y$ if, for all morphisms $h: W \to X$, the homomorphism
$$\bB(W \overset h \to\longrightarrow X) \overset {\bullet \theta} \to\longrightarrow 
\bB(W \overset {f\circ h} \to\longrightarrow Y)$$
is an isomorphism. If $\bB$ is graded and $\theta \in \bB^d(X \overset f \to\longrightarrow Y)$, then we say that $\theta$ is of codimension $d$ or of dimension $-d$. And if the morphism $f$ has a strong orientation in the bivariant group $\bB(X \overset f \to\longrightarrow Y)$ it is said to be {\it strongly orientable with respect to the bivariant theory $\bB$}.
\enddefinition

All the Grothendieck transformations constructed in [FM] extend simple functors on the corresponding contravariant functors, but there are few uniqueness theorems available (see [FM, 10.9]). However, the appeal to the strong orientation certainly gives us a uniqueness theorem. Indeed, suppose that we have a Grothendieck transformation between two bivariant theories 
 $$\ga : \bB \to \bB'.$$ 
If for an object $Y$ there is an object $Z$ such that there exists a morphism $g: Y \to Z$ and a certain bivariant element $e_g \in \bB(Y \overset g \to \longrightarrow Z)$ such that $\ga_g (e_g) \in \bB'(Y \overset \ g \to \longrightarrow Z)$ is a strong orientation for the morphism $ g$ and furthermore for morphisms whose target object is $Z$ the Grothendieck transformations $\ga : \bB(? \to Z) \to \bB'(? \to Z)$ are uniquely determined, then one can see that restricted to morphisms $f: X \to Y$ with the target object being $Y$ the Grothendieck transformation is uniquely determined. Indeed, for any element $\alp \in \bB(X \overset f \to \longrightarrow Y)$ we have
$$\ga_{g \circ f} (\alp \bullet e_g) = \ga_f (\alp) {\bullet} \ga_g (e_g). $$
Since $ {\bullet} \ga_g (e_g) : \bB'(X \to Y) \to \bB'(X \to Z)$ is an isomorphism and $\ga_{g \circ f} : \bB(X \overset {g \circ f} \to \longrightarrow Z) \to \bB'(X \overset {g \circ f} \to \longrightarrow Z)$ is uniquely 
determined, it follows that $\ga_f (\alp)$ is uniquely described by
$$\ga_f (\alp) = \Bigl  (\bullet \ga_g (e_g) \Bigr )^{-1} \Bigl  (\ga_{g \circ f} (\alp \bullet e_g)\Bigr ).$$

In what follows, the object $Z$ will be a point.

\proclaim {Proposition (3.2)} Let $Y$ be a possibly singular analytic variety such that the morphism $c:Y \to pt$ has a strong orientation $\theta \in \bH(Y \to pt) = H_*(Y)$ which is contained in the image of $c_*: F(Y) \to H_*(Y)$. Then for any morphism $f: X \to Y$ a bivariant Chern class 
$$\ga_f: \bF(X \overset f \to\longrightarrow Y)  \to \bH(X \overset f \to\longrightarrow Y)$$
is uniquely determined.
\endproclaim

\demo {Proof} First we note that for any variety $X$ the homomorphism 
$\ga_{X \to pt} : F(X) = \bF(X \to pt) \to \bH(X \to pt) = H_*(X)$ is nothing but the Chern--Schwartz--MacPherson class homomorphism $c_*: F(X) \to H_*(X)$. It is because $\ga_{X \to pt}$ is a natural transformation satisfying the normalization condition and thus it has to be the 
Chern--Schwartz--MacPherson class $c_*: F(X) \to H_*(X)$ since it is unique (by resolution of singularities).

Since the strong orientation $\theta$ belongs to the image of $c_*: F(X) \to H_*(X)$, there exists a certain constructible function $\sigma_{\theta} \in F(Y)$ (not necessarily uniquely determined) such that
$$ \theta = c_*(\sigma_{\theta}) = \ga_{Y \to pt} (\sigma_{\theta}).$$

Now, for any bivariant constructible function$\alpha \in \bF(X \overset f \to\longrightarrow Y)$ we have
$$\ga_{X \to pt} (\alpha \bullet \sigma_{\theta}) = \ga_f (\alpha) \bullet \ga_{Y \to pt} (\sigma_{\theta}),$$
which means that 
$$c_*( \alpha \bullet \sigma_{\theta}) = \ga_f (\alpha) \bullet \theta.$$
Therefore it follows that 
$$\ga_f (\alpha) = (\bullet \theta)^{-1}\Bigl  (c_*( \alpha \bullet \sigma_{\theta}) \Bigr ),$$
namely, the homomorphism$\ga_f: \bF(X \overset f \to\longrightarrow Y)  \to \bH(X \overset f \to\longrightarrow Y)$ is uniquely determined. \qed
\enddemo
 
\remark {Remark (3.3)} Note that for any constructible function $\be \in F(Y')$ and a base change morphism $g:Y' \to Y$, we have
$$c_*(g^*\alpha \bullet \be') = g^* \bigl (\ga_f(\alpha) \bigr ) \bullet c_*(\be'). $$
The uniqueness of $\ga: \bF \to \bH$, therefore, would follow if we can show 
that $\omega \in \bH(X \overset f \to\longrightarrow Y)$ and 
$g^*\omega \bullet c_*(\be') = 0$ for any $g:Y' \to Y$ and $\be' \in F(Y)$ automatically implies that $\omega = 0$. At the moment we do not know how to show this. (However, for a possible ``solution" to this problem see the last two sections.)
\endremark

Now the question is what kind of variety $Y$ satisfies the condition described in the above proposition, 
i.e., that the morphism $c:Y \to pt$ has a strong orientation $\theta \in \bH(Y \to pt) = H_*(Y)$ which is in the image of $c_*$. For that purpose we appeal to the derived-category aspect of the bivariant homology theory. 

Let $\Cal C$ be the category of locally compact spaces with morphisms having finite cohomological dimension. 
Let $A$ be a fixed Noetherian ring (e.g., like $\bZ$ and $\bQ$) and let $A_X$ denote the constant sheaf $A$ on the space $X$. We set
$$\bSH^i_A(X \overset f \to\longrightarrow Y) := \op {R^iHom}(Rf_!A_X, A_Y).$$
By Verdier duality, we have
$$\op {R^iHom}(Rf_!A_X, A_Y) = \op {R^iHom}(A_X, f^!A_Y) = H^i(X, f^!A_Y).$$
It turns out that this becomes a bivariant theory [FM, 7.3.1], which shall be called {\it the sheaf-theoretic bivariant homology theory} with the coefficients in $A$ and the bivariant homology theory $\bH$ given in the previous section is called the {\it topological bivariant homology theory}. And furthermore we have the following theorem:

\proclaim {Theorem (3.4)} ([FM, 7.3.4, p.86]) On the subcategory of $\Cal C$ consisting of spaces embeddable in Euclidean spaces, the above sheaf-theoretic bivariant theory is isomorphic to the topological bivariant homology theory with the coefficient in $A$:
$$\bSH^i_A(X \overset f \to\longrightarrow Y) \cong  \bH^i(X \overset f \to\longrightarrow Y) \otimes A.$$
\endproclaim

Since we have the following isomorphisms
$$\bSH^{-n}_A(X \overset f \to\longrightarrow Y) = \op {R^{ -n}Hom}(A_X, f^!A_Y) = \op {R^0Hom}(A_X[n], f^!A_Y),$$
the above theorem simply means that the (topological) bivariant homology classes $c \in \bH^{-n}(X \overset f \to\longrightarrow Y)\otimes A$ is in one-to-one
correspondence with homomorphisms $ c : A_X[n] \to f^!A_Y$.

As to the problem of strong orientability of a map, we have the following sheaf-theoretic criterion:

\proclaim {Proposition (3.5)} ([FM, 7.3.2, p.85]) The map $f: X \to Y$ has a strong orientation
in $\bH^{-n}_A(X \overset f \to\longrightarrow Y)$ if and only if $f^{!}A_Y$ is quasi-isomorphic to $A_X[n]$.
(Such a map is called homologically normally nonsingular.)
\endproclaim

In our case, we just consider the constant map $c:Y \to pt$ on a complex analytic variety $Y$, which is assumed to be connected in the following calculation. Then we have the following isomorphism:
$$\bSH^{-2 \op {dim_{\bC}Y}}(Y \overset c \to\longrightarrow pt) = \op {R^{-2 dim_{\bC}Y}Hom}(A_Y, c^! A_{pt}) = \op {R^0Hom}(A_Y[2 \op {dim_{\bC}Y}], c^! A_ {pt}).$$
Thus the topological bivariant homology classes $c \in \bH^{-2 \op {dim_{\bC}Y}}(Y \overset c \to \longrightarrow pt)\otimes A$ is in one-to-one
correspondence with homomorphisms $ c : A_Y [2 \op {dim_{\bC}Y}] \to c^! A_{pt}$.
Now by the definition we have
$$\bH^{-2 \op {dim_{\bC}Y}}(Y \overset c \to \longrightarrow pt) \cong H_{2 \op {dim_{\bC}Y}}(Y; A).$$
Since $Y$ is connected, $H_{2 \op {dim_{\bC}Y}}(Y; A)$ is generated by the fundamental class $[Y]$, i.e.,
$$H_{2 \op {dim_{\bC}Y}}(Y; A) = A [Y].$$
Thus the fundamental class $[Y]$ corresponds to the canonical homomorphism 
$$A_Y[2 \op {dim_{\bC}Y}] \to c^! A_{pt}$$
which is induced by capping with the fundamental class $[Y]$.
The complex $c^! A_{pt}$ is, by definition, the so-called dualizing complex $D_Y \otimes A$ (or sometimes denoted by $\omega_Y \otimes A$). Since we have the following isomorphism at each stalk 
$$\bold H^{-p}(D_Y \otimes A)_y \cong H_p(Y, Y \setminus y; A),$$
the above canonical homomorphism $A_Y[2 \op {dim_{\bC}Y}] \to D_Y \otimes A$ is quasi-isomorphism if and only if $H_{2 \op {dim_{\bC}Y}}(Y, Y \setminus y;A) = A$ and $H_p(Y, Y \setminus y; A) = 0$ for $p \not = 2 \op {dim_{\bC}Y}$, which means that $Y$ is an oriented {\it $A$-homology manifold}. 

 Applying this argument to all connected components, we therefore get the following corollary:

\proclaim {Corollary (3.6)} Let $Y$ be a complex variety. Then the fundamental class $[Y] \in H_*(Y; A) \cong \bH^{-*}(Y \overset c \to\longrightarrow pt)\otimes A$ is a strong orientation if and only if $Y$ is an oriented $A$-homology manifold.
\endproclaim

Let $Y$ be an oriented $A$-homology manifold and define $c^*(Y) \in H^*(Y)$ by $c_*(\jeden_Y) = c^*(Y) \cap [Y]$. Then the $0$-degree part of $c^*(Y)$ is equal to one, so that the cohomology class $c^*(Y)$ is invertible. In particular, $c^*(Y) \cap [Y]$ is also a strong orientation.
\proclaim {Corollary (3.7)} Let $Y$ be a complex analytic variety which is an oriented $A$-homology manifold. If there exists a bivariant Chern class $\ga :\bF \to \bH$, then for any morphism $f: X \to Y$ the bivariant Chern class $\ga_f: \bF(X \overset f \to \longrightarrow Y)\otimes A \to \bH( X \overset f \to \longrightarrow Y)\otimes A$ is uniquely determined and it is described by
$$\ga_f(\alp) = f^*c^*(Y)^{-1}\cap c_*(\alp).$$
\endproclaim

When $Y$ is nonsingular, we see that the cohomolgy class $c^*(Y)$ is nothing but the total Chern class $c(TY)$ of the tangent bundle $TY$, hence the inverse $c^*(Y)^{-1}$ is the total Segre class $s(TY)$. Therefore this twisted class $f^*c^*(Y)^{-1}\cap c_*(\alp)$ shall also be called the Ginzburg--Chern class of $\alp$ and still denoted by $\ga^{\op {Gin}}(\alp)$.

In the case when $A = \bQ$, examples of rational homology manifolds include surfaces with Kleinian singularities, the moduli space for curves of a given genus, and more generally {\it Satake's $V$-manifolds} or {\it orbifolds}. In particular, a quotient variety of a nonsingular variety by a finite group is a rational homology manifold.

Conversely we ask ourselves whether the above Ginzburg--Chern class  becomes a Grothendieck transformation for morphisms whose target varieties are oriented $A$ - homology manifolds.

\proclaim {Theorem (3.8)} For a morphism of complex analytic varieties $f: X \to Y$ with $Y$ a $A$-homology  manifold, we define $ \overline{\bF} (X \overset f \to\longrightarrow Y)$ to be the set of all constructible functions $\alp \in F(X)$
satisfying the following two conditions ($\sharp$) and ($\flat$) : for any fiber square
$$\CD
X' @> g' >> X \\
@V f' VV @VV f V\\
Y' @> g >> Y, \endCD
$$
with $Y'$ a $A$-homology  manifold 

($\sharp$) the following equality holds for any constructible function $\be' \in F(Y')$:
$$\ga^{\op {Gin}} (g^{\star}\alp \bullet \be') = \ga^{\op {Gin}}(g^{\star}\alp) \bullet \ga^{\op {Gin}}(\be'),$$

($\flat$) 
$$\ga^{\op {Gin}}(g^{\star}\alp) = g^{\star} \ga^{\op {Gin}}(\alp).$$

\noindent Then $\overline {\bF}$ becomes a bivariant theory with the same operations as in $s\bF$ and furthermore the transformation 
$$\ga^{\op {Gin}}: \overline {\bF} \to \bH$$
is well-defined and becomes the unique Grothendieck transformation satisfying that $\ga^{\op {Gin}}$ for morphisms to a point is the Chern--Schwartz--MacPherson class transformation $c_*: F \to H_*$. And also $\overline {\bF}(X \to pt) = F(X)$. 
\endproclaim

The proof of the theorem is the same as in [Y4]. Note that to prove $\overline {\bF}(X \to pt) = F(X)$ 
we need the cross product formula of the Chern--Schwartz--MacPherson class due to Kwieci\'nski [K] (cf. [KY]).

In a similar manner, we can show the following uniqueness theorem of the Baum--Fulton--MacPherson's bivariant Riemann--Roch theorem with values in the topological bivariant rational homology:

\proclaim {Theorem (3.9)} The Grothendieck transformation
$$\tau : \bK_{\op {alg}} \to \bH_{\bQ}$$
constructed in [FM, Part II] is unique on morphisms whose target varieties are rational homology manifolds.
\endproclaim

\remark {Remark (3.10)} For a more general treatment and a further developement on the above construction of $\overline {\bF}$, see [Sch2].
\endremark

\head \S4 Alexander varieties \endhead

In [KT] S. Kleiman and A. Thorup introduced the notion of {\it $C_{\bQ}$-orthocyclic schemes} and in [Vi1, Vi2] A. Vistoli call this scheme {\it Alexander scheme}. In this section we consider a ground field $k$ of characteristic zero. Thus by H. Hironaka [H] the resolution of singularities is always available in characteristic zero, which is an important fact in what follows.

Let us consider the {\it Fulton--MacPherson operational bivariant Chow group theory} 

\noindent $A(X \to Y) \otimes \bQ$  [FM], which is sometimes simply called the {\it bivariant intersection theory 
with rational coefficients } [F]. 

Now, since the resolution of singularities is always available in characteristic zero, it follows from [Vi1, (1.1) Definition and remarks after it] and [Vi1, (2.2) Definition] that we can simply define an Alexander scheme as follows:

\definition {Definition (4.1)} An algebraic variety $Y$ over $k$ is said to {\it satisfy Alexander duality} 
or it is called an {\it Alexander variety} if the following evaluation homomorphism is an isomorphism:
$$\op {ev}_Y : A(X \to Y) \to A(X)  \qquad \text {defined by} \qquad \op {ev}_Y(\alp) := \alp ([Y])$$
(Note that the value $\alp ([Y])$ is sometimed denoted by $c \cap [Y]$.)
\enddefinition

In fact, the above definition of Alexander variety is nothing but saying that {\it the fundamental class
$[Y]$ is a strong orientation for the morphism $Y \to pt$ to a point $pt$ with respect to the operational bivariant Chow theory}. 
Indeed, the fundamental class 
$[Y]$ induces the canonical operational bivariant class $c_{[Y]} \in A(Y \to pt)$, which is defined as follows:
For any morphism $h:T \to pt$ the homomorphism
$$c_{[Y]} (h) : A(T) \to A(T \times Y)$$
is defined by
$$c_{[Y]} (h)(\delta) := \delta \times [Y].$$
Now for each morphism $X \to Y$ we have the bivariant product
$$A(X \to Y) \otimes A(Y \to pt) \to A(X \to pt),$$
which induces the homomorphism
$$\bullet c_{[Y]} : A(X \to Y) \to A(X \to pt).$$

Here we recall the following [F, Proposition 17.3.1]:
\proclaim {Proposition (4.2)} The following evaluation homomorphism is an isomorphism:
$$\phi :A(X \to pt) \to A(X) \qquad \text {defined by} \qquad \phi(c) := c([pt]) \in A(X).$$
\endproclaim

Now we can see that this proposition implies that the above $\op {ev}_Y : A(X \to Y) \to A(X)$ is an isomorphism if and only if the homomorphism $ A(X \to Y) \overset \bullet c_{[Y]} \to \longrightarrow A(X \to pt)$ is an isomorphism.
Thus it follows that an algebraic variety $Y$ over $k$ is an Alexander variety if $[Y]$ is a strong orientation 
for the operational bivariant Chow group theory with rational coefficients $A \otimes \bQ$.
So, an Alexander variety corresponds to an algebraic version of the notion of an oriented $\bQ$-homology manifold. In particular, a quotient variety of a nonsingular variety by a finite group is an Alexander variety [Vi1, Prop. 2.11 (i)]. 

In this way we get the algebraic counterparts of Corollary (3.7) and Theorem (3.9). In fact, for the operational bivariant Chow group theory we even get the following general result on uniqueness:
 
\proclaim {Theorem (4.3)} If there exists a bivariant Chern class $\ga : \bF \to \bA$ with values in the bivariant intersection theory $\bA$, then it is unique.
\endproclaim

\demo {Proof} Consider a morphism $f:X \to Y$. There exists a ``smooth" envelope 
$$\pi : \widetilde Y \to Y.$$
Namely $\widetilde Y$ is nonsingular and $\pi : \widetilde Y \to Y$ is an envelope (see [F] and [Ki, Remark 3.2]). Then for the following fiber square
$$\CD
\widetilde X @> {\widetilde \pi} >> X \\
@V {\widetilde f} VV @VV f V\\
\widetilde Y @> \pi >> Y, \endCD
$$
the pullback
$$\pi ^* : \bA(X \overset f \to \longrightarrow Y) \to \bA(\widetilde X \overset {\widetilde f} \to \longrightarrow \widetilde Y)$$
is injective (see [Ki, Lemma 2.1]). So, suppose that there exists a bivariant intersection Chern class $\ga : \bF \to \bA$ and consider the bivariant intersection Chern class $\ga_f (\alp) \in \bF(\bA(X \overset f \to \longrightarrow Y)$ for $\alp \in \bF(X \overset f \to \longrightarrow Y).$ Then for the above smooth envelope
$\pi: \widetilde Y \to Y$ the pullback $\pi^* \ga_f (\alp) = \ga_{\widetilde f}(\pi^* \alp)$ is uniquely determined because the base variety $\widetilde Y$ is nonsingular. Therefore the injectivity of the pullback homomorphism $\pi ^* : \bA(X \overset f \to \longrightarrow Y) \to \bA(\widetilde X \overset {\widetilde f} \to \longrightarrow \widetilde Y)$ implies that the bivariant class $\ga_f(\alp)$ is also uniquely determined. \qed
\enddemo

In the same way we can show the following uniqueness theorem of the bivariant (intersection) Riemann--Roch:

\proclaim {Theorem (4.4)} Fulton--MacPherson's bivariant (intersection) Riemann--Roch transformation $\tau : \bK \to \bA_{\bQ}$ is unique.
\endproclaim

\remark {Remark (4.5)} The operational bivariant Chow group theory treated in [EY1] is not the same as 
Fulton--MacPherson's bivariant intersection theory, but it is a bit coarser in the sense that the compatibility with flat pullback is not required. However, just when it comes to the uniqueness of such a Grothendieck transformation, the compatibility with flat pullback does not get involved at all and the above Theorem (4.3) is already shown in [EY1] without using Chow envelopes. In the same way, Theorem (4.4) can be also proved without using Chow envelopes.
\endremark

\head \S5 Pseudo-bivariant classes and Grothendieck transformations \endhead
In this section we show that in a sense there does exist a unique bivariant Chern class specializing to
the original Chern--Schwartz--MacPherson class $c_*: F \to H_*$ for morphisms to a point.
Even if we are mainly interested in bivariant Chern classes in this paper, the following constructions 
apply in the same way to a much more general bivariant context. For some additional interesting examples, 
we refer to [FM] and [Sch2].

In this section we consider a natural transformation $c_*: F(X)= \bF(X \to pt) \to \bH(X \to {pt}) = 
H_*(X)$ of two covariant theories associated to the bivariant theories $\bF$ and
$\bH$. We assume that $c_*$ commutes with the induced exterior products and that $c_*$ maps
 the unit $\jeden_{pt} \in F(pt)$ to the unit $1_{{pt}} \in H_*({pt})$.

As our guiding example, for $\bF$ 
we take the simple bivariant theory of constructible functions, for $\bH$ we take the bivariant homology theory and for $c_*$ we take the Chern--Schwartz--MacPherson transformation.

Let us now start with the general context (where the notions in bracket refer to this Chern class example).

\definition {Definition (5.1)} We define
$$\widetilde {\bF}(X \overset f \to\longrightarrow Y)$$
to be the set consisting of all $\alp \in \bF(X \overset f \to\longrightarrow Y)$ satisfying the following condition: there exists a bivariant class $B_{\alp} \in \bH(X \overset f \to\longrightarrow Y)$ such that for any base change $g: Y' \to Y$ (without any requirement) of an independent square
$$\CD
X'@> {g'} >> X \\
@V {f'} VV @VV f V\\
Y' @> g>> Y, \endCD
$$
and for any $\be' \in F(Y')$ the following equality holds:
$$c_*(g^*\alp \bullet \be') = {g}^*B_{\alp} \bullet c_*(\be').$$
\enddefinition

The above bivariant class $B_{\alp}$ should ideally be the unique bivariant (Chern) class of $\alp$. However, so far we still do not know if it is the case or not. So, provisionally we call $B_{\alp}$ {\it a pseudo-bivariant (Chern) class of $\alp$}.

\proclaim {Theorem (5.2)} The above $\widetilde {\bF}$ is a bivariant theory. Furthermore $\widetilde {\bF}(X \to pt) = F(X)$.
\endproclaim
\demo {Proof} The proof of the second statement follows from the fact that $c_*:F \to H_*$ commutes with the cross product $\times$. Indeed, consider the fiber square 
$$\CD
X'@> {g'} >> X \\
@V {f'} VV @VV {} V\\
Y' @> g>> pt.\endCD
$$
Then, for any $\alp \in F(X)$ and for any $\be' \in F(Y')$  the following holds:
$$
\align
c_*(g^*\alp \bullet \be') & = c_*(\alp \times \be') \\
& = c_*(\alp) \times c_*(\be') \\
& = g ^* c_*(\alp)) \bullet c_*(\be').
\endalign
$$
Hence, we can take $c_*(\alp)$ for $B_{\alp}$.

For the first statement it suffices to show that the three operations are well-defined,  i.e., that the above imposed condition is preserved under the three operations.

(i) the product operation is well-defined: Let $\alp \in \widetilde {\bF}(X \overset f \to\longrightarrow Y)$ and $\be \in \widetilde {\bF}(Y \overset g \to\longrightarrow Z)$. For a base change $h:Z' \to Z$ consider the following independent squares:
$$\CD
X'@> {h''} >> X \\
@V {f'} VV @VV f V\\
Y' @> {h'}>> Y\\ 
@V {g'} VV @VV g V\\
Z' @> {h}>> Z. \endCD \tag 5.2.1
$$
$\alp \in \widetilde {\bF}(X \overset f \to\longrightarrow Y)$ implies that there exists a class $B_{\alp} \in \bH(X \overset f \to\longrightarrow Y)$ such that for any constructible function $\be' \in F(Y')$
$$c_*({h'}^*\alp \bullet \be') = {h'}^*B_{\alp} \bullet c_*(\be') \tag 5.2.2$$
and $\be \in \widetilde {\bF}(Y \overset g \to\longrightarrow Z)$ implies that there exists a class $B_{\be} \in \bH(Y \overset g \to\longrightarrow Z)$ such that for any constructible function $\de' \in F(Z')$
$$c_*(h^*\be \bullet \de') = h^* B_{\be} \bullet c_*(\de').$$
Since $h^*\be \bullet \de' \in F(Y')$, it follows from (5.2.2) that we get
$$
\align
c_* \Bigl ({h'}^*\alp \bullet (h^*\be \bullet \de') \Bigr ) & = {h'}^* B_{\alp} \bullet c_*(h^*\be \bullet \de') \\
& = {h'}^*B_{\alp} \bullet \Bigl  (h^*B_{\be} \bullet c_*(\de') \Bigr ). \endalign
$$
Since ${h'}^*\alp \bullet (h^*\be \bullet \de') = ({h'}^*\alp \bullet h^*\be) \bullet \de' = h^*(\alp \bullet \be) \bullet \de'$, we get that
$$
\align
c_*\Bigl (h^*(\alp \bullet \be) \bullet \de' \Bigr )  & =  \Bigl ({h'}^*B_{\alp} \bullet h^* B_{\be} \Bigr ) \bullet c_*(\de') \\
& = h ^* (B_{\alp} \bullet B_{\be}) \bullet c_*(\de '). \endalign
$$
Therefore we get that $\alp \bullet \be \in  \widetilde {\bF}(X \overset {gf} \to\longrightarrow Z)$, in other words
$B_{\alp} \bullet B_{\be}$ is a pseudo-bivariant (Chern) class of $\alp \bullet \be$.

(ii) the pushforward operation is well-defined: Consider the independent squares (5.2.1) above with $f$ confined and let $\alp \in \widetilde {\bF}(X \overset {gf} \to\longrightarrow Z)$. Then there exists a bivariant class $B_{\alp} \in \bH(X \overset {g f} \to\longrightarrow Z)$ such that for any constructible function $\de' \in F(Z')$
$$c_*(h^*\alp \bullet \de') = h^*B_{\alp} \bullet c_*(\de').$$
From which we get that
$$ {f'}_* \Bigl  (c_*(h^*\alp \bullet \de') \Bigr ) = 
{f'}_* \Bigl (h^*B_{\alp} \bullet c_*(\de') \Bigr ) .$$
From the naturality of $c_*$  we get that 
$$
\align
{f'}_* \Bigl  (c_*(h^*\alp \bullet \de') \Bigr ) &= c_* \Bigl  (f'_*(h^*\alp \bullet \de') \Bigr ) \\
&= c_* \Bigl  ((f'_*h^*\alp) \bullet \de' \Bigr ) \\
&= c_* \Bigl  (h^*(f_*\alp) \bullet \de' \Bigr ). \endalign
$$
Since $ {f'}_* \Bigl  (h^*B_{\alp} \bullet c_*(\de') \Bigr ) =  \Bigl  ({f'}_* h^* B_{\alp}\Bigr )  \bullet c_*(\de')$, we get that
$$
\align
c_* \Bigl  (h^*(f_*\alp) \bullet \de' \Bigr ) & =  \Bigl  ({f'}_* h^* B_{\alp}\Bigr )  \bullet c_*(\de') \\
& = h^* (f_*B_{\alp}) \bullet c_*(\de'). \endalign
$$
Therefore we get that $f_*\alp \in  \widetilde {\bF}(Y \overset {g} \to\longrightarrow Z)$ with $ f_*B_{\alp}$ being a pseudo-bivariant (Chern) class of $f_*\alp$.

(iii) the pullback operation is well-defined: Consider the independent squares
$$\CD
X'' @> {h'} >> X' @> {g'} >> X \\
@V {f''} VV @VV {f'} V@VV f V\\
Y'' @> h >> Y'@> g >> Y. \endCD \tag 5.2.3
$$
Let  $\alp \in \widetilde {\bF}(X \overset {f} \to\longrightarrow Y)$. Then, for the base change $gh: Y'' \to Y$, there exists a class $B_{\alp} \in \bH(X \overset f \to\longrightarrow Y)$ such that for any constructible function $\be'' \in F(Y'')$
$$c_* \Bigl ((gh)^*\alp \bullet \be'' \Bigr ) = (g h)^*B_{\alp} \bullet c_*(\be''),$$
which can be written as
$$c_* \Bigl (h^*(g^*\alp) \bullet \be'' \Bigr ) = h^*( g^*B_{\alp)}) \bullet c_*(\be'').$$
Therefore we get that $g^*\alp \in \widetilde {\bF}(X' \overset {f'} \to\longrightarrow Y')$ with $ g ^* B_{\alp}$ being a pseudo-bivariant (Chern) class of $g^*\alp$.
\qed
\enddemo

Let us illustrate this construction in our guiding case of the Chern -- Schwartz -- MacPherson transformation.

\noindent {\bf Example (5.3):} Let $f:X \to Y$ be a smooth morphism of possibly singular varieties. Then we have
$$\jeden_X \in \widetilde {\bF}(X \overset f \to\longrightarrow Y)$$
with $c^*(T_f) \bullet [f]$ being a pseudo-bivariant Chern class of $\jeden_X$. 
Here $T_f$ is the vector bundle of tangent spaces of fibers of $f$ and $[f] \in \bH(X \overset f \to\longrightarrow Y)$ is the canonical orientation of the smooth morphism $f$.

Indeed, consider the following fiber squares:
$$
\CD
X' @> {} >> X \\
@V {id} VV @VV {id} V\\
X' @> {} >> X \\ 
@V {f'} VV @VV {f} V\\
Y' @> {g} >> Y.\\
\endCD
$$
Then for any $\be' \in F(Y')$ we have
$$
\align
c_*(g^*\jeden_X\bullet \be') &= c_*(\jeden_{X'} \bullet \be') \\
& = c_*({f'}^*\be') \\
& = c(T_{f'}) \cap {f'}^!c_*(\be') \\
& = c(T_{f'}) \bullet [f'] \bullet c_*(\be') \\
& = g^*c(T_f) \bullet g^*[f] \bullet c_*(\be') \\
& = g^*\Bigl  (c(T_f)\bullet [f] \Bigr ) \bullet c_*(\be'). \endalign
$$
The third equality above follows from the so-called Verdier--Riemann--Roch theorem for Chern class ([FM], [Sch1], [Y1]).
\remark {Remarks (5.4)} (1) Let $f:X \to Y$ be a cellular morphism of complex varieties embeddable into complex manifolds. Then $\bF(X \overset f \to\longrightarrow Y) \subset \widetilde {\bF}(X \overset f \to\longrightarrow Y)$ by Theorem (2.1) (if we consider only such embeddable varieties).

(2) Let $f:X \to Y$ be any morphism with $Y$ being nonsingular. Then we will show in a future paper that $\bF(X \overset f \to\longrightarrow Y) \subset \widetilde {\bF}(X \overset f \to\longrightarrow Y)$ in this case also, based on the {\it generalized Verdier--Riemann--Roch theorem for Chern class} [Sch1].
\endremark

Let us now come back to the general bivariant context.
Suppose that $\alp \in \widetilde {\bF}(X \to Y)$. Then there exists a class $B_{\alp} \in \bH(X \to Y)$
such that for any $g: Y' \to Y$ and for any $\be' \in F(Y')$ the following holds
$$c_*(g^*\alp \bullet \be) = g^*B_{\alp} \bullet c_*(\be').$$
Here we note that the bivariant class $B_{\alp}$ might not be uniquely determined. 
And as we saw in the previous sections, in some cases it is surely uniquely determined, 
but in general we still do not know if it is the case. So, to remedy this unpleasant possible non-uniqueness of the bivariant class $B_{\alp}$, we set
$$
\align
& \bPH(X \overset f \to\longrightarrow Y) := \\
& \Bigl  \{ B \in \bH(X \overset f \to\longrightarrow Y) | \text {$B$ is a pseudo-bivariant (Chern) class of
some $\alp \in \widetilde {\bF}(X \overset f \to\longrightarrow Y)$ } \Bigr \}\endalign
$$
be the set of all pseudo-bivariant (Chern) classes for the morphism $ f: X \to Y$.
It is clear that $\bPH$ is a bivariant subtheory of $\bH$, i.e, it is stable under the three bivariant operations, which can be seen from the proof of Theorem (5.2).

Then we define
$$\widetilde {\bH}(X \overset f \to\longrightarrow Y) 
:= \bPH(X \overset f \to\longrightarrow Y) / \sim$$
where the relation $\sim$ is defined by
$$B \sim B' \Longleftrightarrow g^*B \bullet c_*(\be') = g^*B' 
\bullet c_*(\be')$$
for all independent squares with $g:Y' \to Y$ and all $\be' \in F(Y')$.

We note that if there exists $e_Y \in F(Y)$ such that $c_*(e_Y)$ is a strong orientation, then the equivalence relation $\sim$ on $\bPH(X \to Y)$ is the identity relation. In particular, this is the case for $Y = pt$ the final object, since $c_*$ preserves the corresponding units by assumption.

Certainly the relation $\sim$ is an equivalence relation. In other words, with this identification we want to make possibly many pseudo-bivariant (Chern) classes into one unique bivariant (Chern) class. Indeed we have

\proclaim {Theorem (5.5)} $\widetilde {\bH}(X \overset f \to\longrightarrow Y)$ is an Abelian group and $\widetilde {\bH}$ is a bivariant theory with the canonical operations induced from those of $\bH$. Furthermore we have
$$\widetilde {\bH}(X \to {pt}) = \op {Image}(c_*:F(X) \to H_*(X)).$$
\endproclaim
\demo {Proof} Let us denote the equivalence class of $B \in \bPH(X \to Y)$ by $[B]$. 

It is easy to see that $\widetilde {\bH}(X \overset f \to\longrightarrow Y)$ is an Abelian group, because the definition of the sum $[B_1] + [B_2]:= [B_1 + B_2]$ is well-defined. Indeed, suppose that $B_1 \sim B'_1$ and $B_2 \sim B'_2$, i.e.,
$$ g^*B_1 \bullet c_*(\be') = g^*B_1' \bullet c_*(\be'), 
\quad g^*B_2 \bullet c_*(\be') = g^*B_2' \bullet c_*(\be')$$
for any $g:Y' \to Y$ and for any $\be' \in F(Y')$. Then it is obvious that 
$$ g^*(B_1+ B_2) \bullet c_*(\be') = g^*(B'_1+ B'_2) \bullet c_*(\be'). $$
Thus we get that $(B_1 +B_2) \sim (B'_1 + B'_2).$

Now let us see that $\widetilde {\bH}$ is a bivariant theory with the canonical operations induced from those of $\bH$. Define the product
$$ \widetilde {\bullet} : \widetilde {\bH}(X \overset f \to\longrightarrow Y) 
\otimes \widetilde {\bH}(Y \overset g \to\longrightarrow Z) \to 
\widetilde {\bH}(X \overset {g f} \to\longrightarrow Z)$$
by
$$[B_1] \widetilde {\bullet} [B_2] := [B_1 \bullet B_2].$$
The well-definedness of the product, i.e, $B_1 \sim B'_1$ and $B_2 \sim B'_2$ imply $B_1 \bullet B_2 \sim B'_1 \bullet B'_2$, can be seen as follows (cf. the proof of Theorem (5.2)(i)):
$$ 
\align
h^*(B_1 \bullet B_2) \bullet c_*(\de') 
&= \Bigl  ({h'}^*B_1 \bullet h^*B_2 \Bigr ) \bullet c_*(\de') \\
& =   {h'}^*B_1 \bullet \Bigl  ( h^*B_2 \bullet c_*(\de') \Bigr ) \\
& =  {h'}^*B_1 \bullet \Bigl  ( h^*B'_2 \bullet c_*(\de') \Bigr ) \\
& =  {h'}^*B'_1 \bullet \Bigl  ( h^*B'_2 \bullet c_*(\de') \Bigr ) \\
& = \Bigl  ({h'}^*B'_1 \bullet h^*B'_2 \Bigr ) \bullet c_*(\de')  \\
& = h^*(B'_1 \bullet B'_2) \bullet c_*(\de')  . \endalign
$$
In the fourth equality above we use the equivalenc relation $B_1 \sim B'_1$ since $B'_2$ is a pseudo-bivariant class so that
$ h^*B'_2 \bullet c_*(\de') = c_*(\de'')$ for some $\de'' \in F(Y')$.

The well-definedness of the pushforward (for $f$ confined)
$$ f_{\bigstar}:\widetilde {\bH} (X \overset {g f} \to\longrightarrow Z) \to \widetilde {\bH}(Y \overset {g} \to\longrightarrow Z) \quad \text {defined by} \quad  f_{\bigstar}[B] := [f_*B],$$
 i.e., that $B \sim B'$ implies $ f_*B \sim f_* B'$, can be seen as follows. We use the diagram (5.2.1) above.$B \sim B'$ means that
$ h^*B \bullet c_*(\de') = h^*B' \bullet c_*(\de') $ for any $h: Z' \to Z$ and for any $\de' \in F(Z')$. Then we have
$$
\align
h^*(f_*B) \bullet c_*(\de') & = {f'}_*(h^*B) \bullet c_*(\de') \\
& =  {f'}_* \Bigl  (h^*B \bullet c_*(\de') \Bigr ) \\
& =  {f'}_* \Bigl  (h^*B' \bullet c_*(\de') \Bigr ) \\
& = {f'}_*(h^*B') \bullet c_*(\de') \\
& = h^*(f_*B') \bullet c_*(\de'). \endalign
$$
Therefore we get that $ f_*B \sim f_*B'$.

The pullback
$$ g^{\bigstar}:\widetilde {\bH} (X \overset {f} \to\longrightarrow Y) \to \widetilde {\bH} (X' \overset {f'} \to\longrightarrow Y)$$
is defined by
$$ g^{\bigstar} [B] := [g^*B].$$
Its well-definednes can be seen as follows. Consider the diagram (5.2.3). We want to show that $B \sim B'$ implies $ g^*B \sim g^*B'$. Let $h: Y'' \to Y'$ be a base change and let $\be'' \in F(Y'')$.
Then, since $gh: Y'' \to Y$ is a base change of $Y$ and $B \sim B'$, we have that $(g h)^*B \bullet c_*(\be'') = (g h)^*B' \bullet c_*(\be'').$ Hence
$$
\align
h^*(g^*B) \bullet c_*(\be'') &= (g h)^*B \bullet c_*(\be'') \\
& = (g h)^*B' \bullet c_*(\be'') \\
& = h^*(g^*B') \bullet c_*(\be'').\endalign
$$
Therefore we get that $ g^*B \sim g^*B'$.

It is left for the reader to show that $\widetilde {\bH}(X \to {pt}) = \op {Image}(c_*:F(X) \to H_*(X))$. More precisely, using the cross product formula of $c_*$ one can show that 
$c_*: \widetilde {\bF}(X \to pt) = F(X \to pt) \to \bPH(X \to {pt})$ is surjective and 
we already pointed out before that the equivalence relation $\sim$ on $\bPH(X \to {pt})$ is the identity relation.
 \qed
\enddemo

\remark {Remark (5.6)} In our previous paper [Y3] we posed the problem of whether or not there is a reasonable bivariant homology theory so that the following ``quotient"
$$\frac {c_*(\alp)}{c_*(\jeden_Y)}$$
is well-defined. The above theory $\widetilde {\bH}$ is in a sense a positive answer to this problem.
\endremark

Now it is easy to see the following
\proclaim {Theorem (5.7)} There exists a unique Grothendieck transformation
$$\widetilde {\ga} : \widetilde {\bF} \to \widetilde {\bH}$$
whose associated covariant transformation is $c_*:F \to \op {Im}(c_*)$, where $\op {Im}(c_*)(X) := \op {Image}\Bigl (c_*: F(X) \to H_*(X) \Bigr )$.
\endproclaim

\demo {Proof} Clearly the existence follows from defining $\widetilde {\ga} : \widetilde {\bF}(X \to Y)  \to \widetilde {\bH}(X \to Y)$ by
$$\widetilde {\ga} (\alp) := [B_{\alp}]$$
where $B_{\alp}$ is a pseudo-bivariant (Chern) class of $\alp$.
This is well-defined. Indeed, suppose that $B_{\alp}$ and $B'_{\alp}$ are pseudo-bivariant (Chern) classes of $\alp$, which means that for any independent square with base change $g:Y' \to Y$ and for any $\be' \in F(Y')$
$$ g^* B'_{\alp} \bullet c_*(\be') = c_*(g^*\alp \bullet \be') = g^* B_{\alp} \bullet c_*(\be').$$ 
Thus clearly we get that 
$$ g^* B_{\alp} \bullet c_*(\be') = g^* B'_{\alp} \bullet c_*(\be')$$
for any independent square with base change $g:Y' \to Y$ and for any $\be' \in F(Y')$, i.e., $B_{\alp} \sim B'_{\alp}$. Moreover, it follows from Theorem (5.5) and the proof of Theorem (5.2) that $\widetilde {\ga}$ is a Grothendieck transformation.

The uniqueness follows from the construction. Indeed, let $\widetilde {\ga}': \widetilde {\bF} \to \widetilde {\bH}$ be a
Grothendieck transformation whose associated covariant transformation is the given $c_*$. Then, for $\alp \in \widetilde {\bF}(X \to Y)$ we have
$$\widetilde {\ga}'(\alp) = [B']$$
for some pseudo-bivariant (Chern) class $B' \in \bPH(X \to Y)$ (of some $\alp' \in \widetilde {\bF}(X \to Y)$). Using the fact that the equivalence relation $\sim$ on $\bPH(X \to Y)$ is the identity relation when $Y = pt$ is a final object,  for any base change $g: Y' \to Y$ and for any $\be' \in F(Y')$ we have
$$
\align
g^*B_{\alp} \bullet c_*(\be') & = c_*(g^*\alp \bullet \be') \\
& = [c_*(g^* \bullet \be')] \\
& = \widetilde {\ga}'(g^*\alp \bullet \be') \\
& = \widetilde {\ga}'(g^*\alp) \widetilde \bullet \widetilde {\ga}'(\be') \\
& = g^{\bigstar}\widetilde {\ga}'(\alp) \widetilde \bullet  [c_*(\be')] \\
& = g^{\bigstar}[B'] \widetilde \bullet  [c_*(\be')] \\
& = [g^*B'] \widetilde \bullet  [c_*(\be')] \\
& = [g^*B' \bullet  c_*(\be')] \\
& = g^*B'  \bullet  c_*(\be'). \endalign
$$
Hence we get that $B_{\alp} \sim B'$, which means that $\widetilde {\ga}(\alp) = \widetilde {\ga}'(\alp)$.
Thus the uniquenss follows.
\qed
\enddemo

\remark {Remark (5.8)} If $c_*:F \to H$ is the Chern--Schwartz--MacPherson class in the above construction, the above theorem can be put it as follows: There exists a unique Grothendieck transformation
$$\widetilde {\ga} : \widetilde {\bF} \to \widetilde {\bH}$$
satisfying the normalization condition that for a nonsingular variety $X$
$$\widetilde {\ga} (\jeden_X) = c(TX) \cap [X].$$
\endremark

\remark {Remark (5.9)} The above construction gives in some sense a positive answer to the general uniqueness problem of Grothendieck transformations [FM, \S 10.9 Uniqueness questions]. 
Suppose that there is a Grothendieck transformation $\ga: \bF' \to \bH$ for a bivariant subtheory $\bF'$ of $\bF$ such that its associated covariant transformation is the given $c_*$. Then it is clear from the definition of pseudo-bivariant classes that the image of $\ga$ is containd in the subtheory $\bPH$ of $\bH$. The Grothendieck transformation $\ga: \bF' \to \bH$ restricted to $\bPH$ on the  target side shall also be denoted by $\ga: \bF' \to \bPH$. Let $\iota_{\bF}: \bF' \to \widetilde{\bF}$ be the inclusion transformation and let $q:\bPH \to  \widetilde {\bH}= \bPH/{\sim}$ be the quotient transformation. Then the
following diagram commutes:
$$\CD
\bF'@> {\ga} >> \bPH \\
@V {\iota_{\bF}} VV @VV q V\\
\widetilde {\bF}@> {\widetilde {\ga}}>> \widetilde {\bH}= \bPH/{\sim}. \endCD
$$
So the existence (or uniqueness) of such a Grothendieck transformation
is equivalent to the existence (or uniqueness) of a ``functorial section"
of the quotient transformation $q:\bPH \to \widetilde {\bH}= \bPH/{\sim}$.
This is, clearly, the case if the equivalence relation $\sim$ is the identity relation
on all the groups $\bPH(X \overset {f} \to\longrightarrow Y)$, e.g., as in Theorems (4.3) and (4.4).
If $\iota_{\bH}: \bPH \to \bH$ denotes the inclusion transformation, then  the composite $\iota_{\bH} \circ \tilde{\ga} : \widetilde {\bF} \to \widetilde {\bH} \simeq \bPH  \to  \bH$
is the unique Grothendieck transformation to $\bH$ which extends
the given covariant transformation $c_*$ and such that the subtheory $\widetilde {\bF}$ of $\bF$ is the maximal one.
\endremark   

\remark {Remark (5.10)}
So far we have never used the {\it bivariant projection formula}, i.e., the axiom $A_{123}$ in
the bivariant theory of [FM], so that all our arguments of this paper
apply also to {\it a weak bivariant theory} in the sense of [Sch2].
Similarly, our construction above works mutatis mutandis for
(weak) partial bivariant theories in the sense of [Sch2],
where one associates a bivariant group $\bB(X \overset {f} \to\longrightarrow Y)$ only to a suitable class of 
so-called ``allowable maps". 

For example, let us consider only such maps $f: X \to Y$ that 
$H_*(Y)=\bH(Y \to pt)$ contains a strong orientation in the image of $c_*$
(as in \S3). Then the equivalence relation $\sim$ is the identity relation on $\bPH(X \overset {f} \to\longrightarrow Y)$
and the partial Grothendieck transformation 
$\iota_{\bH} \circ \tilde{\ga} : \widetilde {\bF} \to \widetilde {\bH} \simeq \bPH  \to  \bH$
reduces exactly to the corresponding transformation
constructed in [Sch2, Theorem (2.1)]. In particular in the context
of the Chern--Schwartz--MacPherson transformation $c_*$ with $Y$ being a
smooth (or homology) manifold, we recover the main result of [Y4]
or our Theorem (3.8).
\endremark

\head \S6 A viewpoint from Operational Bivariant Theories \endhead

In this final section we make a few more remarks on our constructions given in the previous section. 

If the independent squares are exactly the Cartesian squares
in our underlying category $\Cal C$, then we can explain our previous constructions 
in the language of the operational bivariant theories. Here we refer to
[FM] for the general context of operational bivariant theories (cf. [EY1], [EY2] and [Y2]) .

Let $B_*$ be a covariant functor (or sometimes called a homology theory) on the category $\Cal C$. Then the {\it associated operational bivariant theory} $\bB^{\op {op}}$ of $B_*$ is defined as follows. For a morphism $f: X \to Y$, an element $c \in \bB^{\op {op}}( X \overset f \to \longrightarrow Y)$ is defined to be a collection of homomorphisms
$$ c(g): B_*(Y') \to B_*(X')$$
for all $g: Y' \to Y$ and the fiber square 
$$\CD
X' @> g' >> X \\
@V f' VV @VV f V \\
Y' @> g >> Y.\endCD 
$$
and these homomorphisms $\{c(g) \}$ are required to be compatible with proper pushforwards in the obvious sense. Then $\bB^{\op {op}}$ can be made into a bivariant theory. Note that for
the definition of the operational pushforward one needs the assumption that 
the independent squares are exactly the Cartesian squares.

If $\Cal C$ has a final object $pt$ and $B_*(pt)$ has a distinguished element $1_{pt}$, then
the homomorphism 
$$\op {ev}: \bB^{\op {op}}(X \to pt) \to B_*(X) \quad \text {defined by } \quad \op {ev}(c) := \bigl (c(\op {id}_{pt}) \bigr )(1_{pt})$$
 is called the {\it evaluation homomorphism} as in \S 4..

Let $\bB$ be a bivariant theory. Then the {\it associated operational bivariant theory} of $\bB$ , denoted also by $\bB^{\op {op}}$, is defined to be the associated operational bivariant theory of the covariant functor $B_*(X) := \bB(X \to pt)$ as above.
Then we have the following canonical Grothendieck transformation
$$\op {op} : \bB \to \bB^{\op {op}}$$
defined by, for each $b \in \bB(X \to Y)$,
$$\op {op}(b) := \Bigl \{ (g^*b)\bullet : B_*(Y') \to B_*(X') | g: Y' \to Y \Bigr  \}.$$ 

For $Y=pt$ we get that $\op {op}(b)$ is the exterior product with $b \in \bB_*(X)$. In particular, $\op {ev} \circ \op {op}$ (of course restricted to morphisms to the final object $pt$)
is the identity so that the covariant transformation $\op {op}: \bB(X \to pt) \to \bB^{\op {op}}(X \to pt)$ is injective.

In our application of operational bivariant theories, the covariant theories $B_*$ which we treat are
$H_*(X) = \bH(X \to pt)$ and $\op {Im}(c_*)(X)$. 

Note that in general $(g^*b) \bullet $ does not preserve the subgroups $\op {Im}(c_*)$.
But by the definition of $\bPH$ this is the case for a pseudo-bivariant 
class $B \in \bPH(X \overset {f} \to\longrightarrow Y)$. Therefore, by restricting $\op {op}: \bH \to \bH^{\op {op}}$ to $\bPH$ on the source side and to $\bigl (\op {Im}(c_*) \bigr )^{\op {op}}$ on the target  side, we get the following Grothendieck transformation
$$\op {op}_{\bPH}: \bPH \to \Bigl (\op {Im}(c_*) \Bigr )^{\op {op}}.$$
Moreover, by the definition of the equivalence relation $\sim$, 
we get that for $B, B' \in \bPH(X \overset {f} \to\longrightarrow Y)$
$$B \sim B' \quad \text {if and only if } \quad \op {op}_{\bPH}(B)= \op {op}_{\bPH}(B').$$
Therefore we get the following injective
Grothendieck transformation 
$$\widetilde {\op {op}}_{\bPH}: \widetilde {\bH}= \bPH/{\sim}  \to  \Bigl (\op {Im}(c_*) \Bigr )^{\op {op}},$$
which thus identifies $\widetilde {\bH}$ with a subtheory of $\bigl (\op {Im}(c_*) \bigr )^{\op {op}}$.

In this way the Grothendieck transformation $\widetilde {\ga} : \widetilde {\bF} \to \widetilde {\bH}$ of Theorem (5.7) induces also a Grothendieck transformation 
$$\widetilde {\ga}^{\op {op}} := \widetilde {\op {op}}_{\bPH} \circ  \widetilde {\ga}: \widetilde {\bF} \to  \Bigl (\op {Im}(c_*) \Bigr )^{\op {op}}$$
such that (restricted to morphisms to the final object) 
$\op {ev} \circ \widetilde {\ga}^{\op {op}}: F \to \op {Im}(c_*)$ is the
given covariant transformation $c_*:F \to \op {Im}(c_*)$.

The uniqueness of such a Grothendieck transformation $\widetilde {\ga}^{\op {op}}$
follows from the fact that the operation $\widetilde {\ga}^{\op {op}}(\alpha) (g)$
fits into the following commutative diagram (cf. [EY1], EY2] and [Y2]):
$$
\CD
F(Y')  @> {g^* \alpha} \bullet >> F(X') \\
@V  {c_*} VV @VV {c_*}V \\
\op {Im}(c_*)(Y') @>  {\widetilde {\ga}^{\op {op}}(\alpha) (g)} >> \op {Im}(c_*)(X') , \endCD 
$$
whose vertical arrows are surjective. Since $\widetilde {\op {op}}_{\bPH}$ is injective and $\widetilde {\ga}^{\op {op}} = \widetilde {\op {op}}_{\bPH} \circ  \widetilde {\ga}$ is unique, it follows that $\widetilde {\ga}$ is also unique. Thus we obtain the uniqueness result of Theorem (5.7). 

\remark {Remark (6.1)} Since $c_*$ commutes by assumption with exterior products,
we can apply Theorem A and Theorem 3.10 of [Y2] to get a Grothendieck
transformation $ \gamma^{\vee} : \bF^{\vee} \to \bigl (\op {Im}(c_*) \bigr )^{\op {op}}$ for some subtheory $\bF^{\vee}$ of $\bF$. It follows from Theorem (5.7) that this subtheory $\bF^{\vee}$ contains our $\widetilde {\bF}$. 
Note that $\widetilde {\bF}$ could be a strictly proper subgroup of $\bF^{\vee}$. The possible difference comes from the fact that
$\widetilde {\ga}^{\op {op}}: \widetilde {\bF} \to  \bigl (\op {Im}(c_*) \bigr )^{\op {op}}$ factorizes through the subtheory $\widetilde {\bH}$ of $\bigl (\op {Im}(c_*) \bigr )^{\op {op}}$. So, in this sense our Grothendieck transformation
$\widetilde {\ga}:\widetilde {\bF} \to \widetilde {\bH}$ interpolates between the (biggest) operational Grothendieck transformation $\ga^{\vee}: \bF^{\vee}\to \bigl (\op {Im}(c_*) \bigr )^{\op {op}}$ (as in [Y2, Theorem 3.10]) and a true Grothendieck transformation
$\ga: \bF' \to \bH$, which may not exist. Even if $\ga$ exists,
it may be not unique. But in any case, it fits into the following commutative
diagrams of Grothendieck transformations:

$$\CD
\bF' @> \ga >> \bPH @> {\op {op}_{\bPH} }>> \Bigl (\op {Im}(c_*) \Bigr )^{\op {op}}\\
@V  {\iota_{\bF}} VV @VV {q}V @VV {\op {id}} V \\
\widetilde {\bF}@>  {\widetilde {\ga}} >> \widetilde {\bH}@>  {\widetilde {\op {op}}_{\bPH} }>> \Bigl (\op {Im}(c_*) \Bigr )^{\op {op}} \endCD 
$$
with  $\widetilde {\op {op}}_{\bPH} \circ  \widetilde {\ga}: \widetilde {\bF} \to  \bigl (\op {Im}(c_*) \bigr )^{\op {op}}$ being the unique Grothendieck transformation,
which factorizes through $\widetilde {\bH}$.
\endremark

\remark {Remark (6.2)} In [EY2] we have constructed a bivariant theory $\widehat {\bF}$ of constructible functions such that there exists a unique bivariant Chern class $\widehat {\ga}: \widehat  {\bF} \to \bA$ satisfying the normalization condition. If we consider the Fulton--MacPherson operational bivariant intersection theory $\bA(X \to Y)$ for $\bH(X \to Y)$ in Definition (5.1), then we can also recover this result
as a special case of our Theorem (5.7), i.e. $\widehat {\bF} = \widetilde {\bF}$ with $\widehat  {\ga}$ being given as the composite $\iota_{\bH} \circ \widetilde {\ga}: \widetilde {\bF} \to \widetilde {\bA} \to \bA$ with $\iota_{\bA}: \widetilde {\bA} \to \bA$ being the inclusion as in Remark (5.9), by appealing to the Chow envelope. The injectivity of $\iota_{\bA}: \widetilde {\bA} \to \bA$, i.e., that the equivalence relation $\sim$ is the identity relation, follows from the fact (which was used in Theorem (4.3)) that for a ``smooth" Chow envelope $\pi: \widetilde Y \to Y$ the bivariant pullback $\pi^*: \bA(X \to Y) \to \bA(\widetilde X \to \widetilde Y)$ is injective. 
Note that the proof in [EY2] is more elementary and the Chow envelopes are not required.
\endremark

\enddocument